\documentclass[12pt]{article}
\usepackage{amsmath,amsfonts}
\usepackage{verbatim}
\usepackage{relsize}
\usepackage{enumerate}
\usepackage[inline]{enumitem}
\usepackage{color}
\normalbaselines
\newcommand\blfootnote[1]{%
  \begingroup
  \renewcommand\thefootnote{}\footnote{#1}%
  \addtocounter{footnote}{-1}%
  \endgroup
}

\begin{document}

\title{{\bf A Gauss-Bonnet-Chern type obstruction for Killing vector fields on Lorentzian
manifolds}\blfootnote{This work has been partially supported by the Spanish MICINN and ERDF project PID2020-116126GB-I00, by the Andalusian and ERDF project A-FQM-494-UGR18, and by the ``Mar\'{\i}a de Maeztu'' Excellence Unit IMAG, reference CEX2020-001105-M, funded by
MCIN-AEI-10.13039-501100011033.} 
}

\author{Alfonso Romero and Miguel S\'anchez
\\ Departamento de Geometr\'{\i}a y Topolog\'{\i}a \\
Facultad de Ciencias, Universidad de Granada \\
 18071-Granada, Spain.\\
E-mail: {\ttfamily aromero@ugr.es, sanchezm@ugr.es}}

\date{  }

\newtheorem{defn}{Definition}

\newtheorem{teo}[defn]{Theorem}
\newtheorem{eje}[defn]{Example}
\newtheorem{counter-exe}[defn]{Counter-example}
\newtheorem{lem}[defn]{Lemma}
\newtheorem{rem}[defn]{Remark}
\newtheorem{cor}[defn]{Corollary}
\newtheorem{pro}[defn]{Proposition}

\font\ddpp=msbm10  at 12 truept
\def\R{\hbox{\ddpp R}}
\def\C{\hbox{\ddpp C}}
\def\L{\hbox{\ddpp L}}
\def\S{\hbox{\ddpp S}}
\def\Z{\hbox{\ddpp Z}}

\newcommand{\D}{{\cal D}}
\newcommand{\M}{{\cal M}}
\newcommand{\Mo}{{\cal M}_0}
\newcommand{\be}{\begin{equation}}
\newcommand{\ee}{\end{equation}}
\newcommand{\la}{\Lambda}
\newcommand{\dimo}{{\bf Proof: }}
\newcommand{\inte}{\int_{0}^{1}}
\newcommand{\gam}{\gamma}
\newcommand{\eps}{\epsilon}
\newcommand{\<}{\langle}
\renewcommand{\>}{\rangle}
\newcommand{\Om}{\Omega^1}
\renewcommand{\(}{\left(}
\renewcommand{\)}{\right)}
\renewcommand{\[}{\left[}
\renewcommand{\]}{\right]}
\newcommand{\om}{\omega}
\newcommand{\me}{\frac{1}{2}}
\newcommand{\Mt}{\widetilde{\M}}
\newcommand{\cat}{{\mathop{\rm cat}\nolimits}}

\hyphenation{Lo-rent-zian}

\maketitle \vspace*{-6mm}
\begin{abstract}
\noindent A new curvature obstruction to the existence of a timelike (resp. causal) Killing or homothetic vector field $X$ on an even-dimensional (odd-dimensional) Lorentzian manifold, in terms of its timelike (resp. null)  sectional curvature is given. As a consequence for the compact case, the well-known Gauss-Bonnet-Chern obstruction to the existence of semi-Riemannian metrics is extended from non-zero constant sectional curvature to non-zero timelike sectional curvature on $X$.
\end{abstract}

\vspace{2mm}

\noindent {\bf 2020 MSC:} Primary 53C50, Secondary 53B30, 53Z05\\
\noindent {\bf Keywords:} Causal and Killing vector fields, timelike sectional curvature, compact Lorentzian manifold. 

\hyphenation{eve-ry}

\hyphenation{di-men-sio-nal}

\thispagestyle{empty}

\section{Introduction}\label{s1}

The presence of a timelike Killing vector field $X$ on an
$m$-dimensional (connected) compact Lorentzian manifold $(M,g)$ has several
interesting consequences on its geometry: 
\begin{enumerate}
\item[1.] It forces the sectional curvature ${\mathcal K}(\pi )$ to be non-positive on any timelike tangent plane $\pi$ containing $X$ in $T_pM$, being $p\in M$ the point where the function $g(X,X)$ attains its maximum value \cite[Th. 3.1(i)]{RS5}.  
\item[2.] It yields further results related to curvature in particular
cases, as when $(M,g)$ is Einstein or it has constant sectional
curvature \cite{RS3,RS5, Ka, Ro}.
\item[3.] It implies geodesic completeness \cite{RS1, RS2}
\end{enumerate}

On the other hand, the case when $(M,g)$ is a compact Lorentzian manifold with constant sectional
curvature $c$ has been especially studied, and specific
obstructions to the values of $c$ can be found:
\begin{enumerate}
\item[{\it i.}] Necessarily, $c\leq 0$. Indeed, such a Lorentzian manifold must be geodesically complete \cite{Ca, Kl}, and, if $c>0$, then $M$ would be a quotient of the Lorentzian pseudosphere $\S^m_1$ (i.e., the $m$-dimensional de Sitter spacetime of constant sectional curvature $c$). But in this case, a result in \cite{CM}, which
states the finiteness of the fundamental group of $M$, would
contradict the compactness of $M$.
\item[{\it ii.}] In the even dimensional case $m=2n$, the Gauss-Bonnet-Chern theorem
implies that no compact Lorentzian manifold $(M,g)$ of constant
sectional curvature $c\neq 0$ exists (see \cite[Prop. 2.8]{Ku2},
for instance).
\item[{\it iii.}] Compact Lorentzian manifolds with constant sectional curvature $c<0$ exist if the
dimension $m$ is odd and, in fact, explicit examples of them with
a timelike Killing vector field  follows from \cite{KR} and the classification in \cite{Ka}.
\end{enumerate}

Thus, it is natural to wonder, beyond the obstructions for  constant sectional curvature, 
\begin{quote}{\it 
Is there a specific relation for even dimensions between the existence of a timelike Killing vector field and the timelike sectional curvature of a compact
Lorentzian manifold? }
\end{quote}

Such an obstruction (including a further extension to the lightlike case) will be obtained in Theorem \ref{t_3} as a consequence of the following result on homothetic vector fields, which has interest in its own right (see Section \ref{s3}). Indeed, it is not restricted to the compact case but it is obviously applicable then.

\begin{teo}\label{tecnica}
Let $(M,g)$ be an $m$-dimensional Lorentzian manifold  and let $X$ be a homothetic vector field on $(M,g)$ such that 
 the function $g(X,X)$ attains a local minimum at $p_{_0} \in
M$ and $X$ is causal at $p_{_0}$ {\rm (}i.e., $g(X_{p_{_0}},X_{p_{_0}})\leq 0, \,X_{p_{_0}}\neq 0${\rm )}. 

{\rm (A)} In case $X_{p_{_0}}$ is timelike and $m$ is even {\rm (}$m=2n, n\geq 1${\rm )},
then $X$ is Killing and there exists a timelike plane
$\pi \subset T_{p_{_0}} M$ such that $X_{p_{_0}} \in \pi$ and its
sectional curvature satisfies
\begin{equation}\label{e_K}
\mathcal{K}(\pi) \geq 0\,.
\end{equation}

{\rm (B)} In case $X_{p_{_0}}$ is lightlike and $m$ is odd {\rm (}$m=2n+1,\, n\geq 1${\rm )},
then there exists a degenerate plane
$\pi \subset T_{p_{_0}} M$ such that $X_{p_{_0}} \in \pi$ and its null sectional curvature with respect to $X_{p_{_0}} $ satisfies 
\begin{equation}\label{e_K2}
\mathcal{K}_{X_{p_{_0}}}(\pi):= 
\frac{g(R(v,X_{p_{_0}})X_{p_{_0}},v)}{g(v,v)}\leq 0\,, 
\end{equation}
where $v$ is any non-lightlike vector in $\pi$.
\end{teo}
\begin{rem}\label{segundo}{\bf (a)} {\rm
{\em About the hypotheses of Theorem \ref{tecnica}}.
Note that the assumptions are imposed only at $p_0$. Thus,  $X$ may be spacelike in some region of $M$ in case (A), see Example \ref{movido1}. In  case (B), $X$ can be even timelike nowhere, see Remark \ref{r_4}(b). 
In case (A), notice that the condition negative sectional curvature on timelike tangent planes does not imply constant sectional curvature (this should happen if the inequality held for all the non-degenerate tangent planes \cite{Ku1}), see \cite[Prop. 4.1]{ ES} or \cite{Sa3} for simple counter-examples).  
In case (B), notice that $\mathcal{K}_{X_{p_{_0}}}(\pi)$  in \eqref{e_K2} is independent of the chosen $v$ and of reversing the sign of $X$ \cite[Prop. 2.3]{Harris}, indeed, its properties are well established in the literature (see \cite[p. 571]{BEE}). On the other hand, inequality \eqref{e_K2} is consistent by continuity with \eqref{e_K} (as the denominator for the computation of $\mathcal{K}(\pi)$ would have different sign for $\mathcal{K}_{X_{p_0}}(\pi)$); in particular, non negative timelike sectional curvature implies \eqref{e_K2}.

{\bf  (b)} {\em About the approach of Theorem \ref{tecnica}}. Our technique is inspired in a result by Berger \cite[p.
169]{P}, which was obtained by studying the classical Hopf
conjecture relating the positiveness of the sectional curvature
and the Euler-Poincar\'e characteristic of a compact even dimensional
Riemannian manifold. In turns, this approach relies on the Bochner technique 
\cite{Wu}, which was suitably adapted for Lorentzian manifolds in \cite{RS3, RS4, Ro}.}
\end{rem}

The following consequence of Theorem \ref{tecnica} extends, on the one hand, the Gauss-Bonnet-Chern obstruction for constant sectional curvature to the non-constant case, when there exists a timelike Killing vector field. Moreover, it gives an obstruction in the compact case to the existence of a Lorentzian metric with negative timelike sectional curvature\footnote{In the non-compact case, the situation is quite different. Indeed, the result \cite[Th. 4]{Ko} shows that any non-compact $n(\geq 2)$ manifold admits a Lorentzian metric with negative timelike sectional curvature. This theorem is obtained by using Gromov's theory of partial differential relations  \cite{Gr} and the $h$-principle on Lorentzian immersions.}.

\vspace{0.5mm}

In the compact semi-Riemannian case any homothetic vector field is Killing (use that the integral of the divergence vanishes) so, only the Killing case is considered at this point.

\begin{teo} \label{t_3}
{\rm (A)} There exists no compact even-dimensional Lorentzian manifold which
admits a timelike Killing vector field $X$ satisfying 
\begin{equation}\label{e_th}
\mathcal{K}(\pi) \neq  0,
\end{equation} for any (timelike) tangent plane $\pi$
which contains $X$, in particular with either everywhere positive or negative
timelike sectional curvature.

{\rm (B)} There exists no compact odd-dimensional Lorentzian manifold 
admitting a causal Killing vector field $X$ which is lightlike at some point and satisfying both \eqref{e_th}
for any timelike tangent plane $\pi$
which contains $X$ and 
\begin{equation}\label{e_th2}
\mathcal{K}_{X}(\pi) \neq  0,
\end{equation} 
for any degenerate tangent plane $\pi$
which contains $X$, 
in particular with either everywhere  positive or negative
timelike\footnote{Notice that no condition on degenerate planes is imposed here, consistently with the last assertion in Remark \ref{segundo}(a).} sectional curvatures. 
\end{teo}

\begin{rem}\label{r_4}{\bf (a)} About the hypotheses of Theorem \ref{t_3} in the timelike case. {\rm 
{\bf (i)}  A timelike Killing vector field on a compact Lorentzian manifold with non negative timelike sectional curvature must be parallel \cite{RS3} (see also \cite[Rem. 5.1]{Sa2}, \cite[Prop. 5.1]{BEM}), and 
{\bf (ii)} A compact manifold admits a Lorentzian metric with a timelike Killing vector field $X$ if and only if it admits an action by $\S^1$ with no fixed points (for the necessary condition, see [RS2, Rem. 5], for the sufficient, see Prop. 5 below); compactness is necessary as $\L^2$ admits such a $X$ but $\R^2$ does not admit a such an action.      

{\bf (b)} {\em About the relation between the cases {\rm (A)} and {\rm (B)} of Theorem \ref{t_3}}. Each timelike Killing vector field $X$ on an $m$-dimensional Lorentzian manifold $(M,g)$ such that $g(X,X)$ attains the global maximum $-c^2$ at some $p_0\in M$, provides a causal Killing vector field,  $\bar X:=X+c\,\partial_\theta$, for the $(m+1)$-dimensional Lorentzian manifold $(M\times \S^1, \bar g:=g+d\theta^2)$, which is lightlike at $p_0$. Observe also that if $c$ is chosen such that  $-c^2$ is the global minimum, one would find an everywhere non-timelike Killing vector field (recall Remark \ref{segundo}(a)).  
}
\end{rem}

In the next section, we will prove Theorems \ref{tecnica} and \ref{t_3}. Additionally, in Section 3, we discuss the extendability of the results to the conformal case and provide further examples and applications, including a relativistic one.
 
\section{Set up and proofs}
Let $(M,g)$ be an $m(\geq 2)$-dimensional (connected) Lorentzian manifold, which we will
consider with signature $(-,+,...,+)$ and let $\nabla$ be the Levi-Civita connection
of $g$. A tangent vector $v$ at
$p\in M$ is said to be {\it timelike} if $g(v,v)<0$, {\it
spacelike} if $g(v,v)>0$ or $v=0$ and {\it lightlike} if
$g(v,v)=0$ and $v\neq 0$.  Consistently, a vector field $X  \in \mathfrak{X}(M)$   on $M$ is 
{\em timelike}  (resp. {\em lightlike}) if so is $X_p$ at each $p\in M$. A tangent (bidimensional) plane $\pi\subset T_pM$ is {\em timelike} (resp. lightlike) if $g_p$ induces a Lorentzian (resp. degenerate) metric on $\pi$. In any of the previous cases, {\em causal} means either  timelike or lightlike. $X$
 is said to be Killing (resp. homothetic) if its
local flows consist in local isometries (resp. homotheties) of $g$. In terms of the
Lie derivative, this is expressed 
$\mathcal{L}_Xg=0$ (resp. $\mathcal{L}_Xg=\lambda g$ for some $\lambda\in \R$). About the existence of a Lorentzian metric with a timelike Killing vector field (recall Remark \ref{r_4}(a), item (ii)), notice:

\begin{pro}\label{existence}
If an $m(\geq 2)$-dimensional manifold $M$ admits an action by the unit circle $\mathbb{S}^1$ without fixed points, then it admits a Lorentzian metric with a timelike Killing vector field.
\end{pro}
\noindent{\it Proof.}  First of all, one can find a Riemannian metric $g_R$ invariant by each transformation $\varphi_t$ induced by the action, namelly $$g_R(u,v):=\int_0^{2\pi}(\varphi^*_tg'_R)(u,v)\,dt,$$ for any $u,v\in T_pM$, $p\in M$, where $g'_R$ is any Riemannian metric, \cite[p. 544]{Wa}.
Therefore, the infinitesimal generator $X$ of $\{\varphi_s\}_{s\in \mathbb{R}}$ is a nowhere zero Killing vector field of the Riemann manifold $(M,g_R)$. 
Now, consider the Lorentzian metric $$g:=g_R-\big(2/g_R(X,X)\big)\,\omega\otimes\omega$$ on $M$, where $\omega$ is the one form on $M$ that is $g_R$-equivalent to $X$. We have $g=g_R$ on $X^{{\perp}_{g_R}}$, $g(v,X)=0$ for any $v\in X^{{\perp}_{g_{_R}}}$ and $g(X,X)=-g_R(X,X)<0$. Thus, $X$ is timelike on $(M,g)$. Finally, note that the vector field $X$ is also Killing on the Lorentzian manifold $(M,g)$. $\hfill\square$	

Now, let  $f:= \frac{1}{2} \,g(X,X)$,  
 and  $A_X$ the (1,1)-tensor field on $M$ defined by
$A_X(v)=-\nabla_vX$, for  any $v\in T_pM$, $p\in M$.  If $X$ is homothetic (or even affine, see \cite[formula (2.3)]{RS5}), then the Hessian of $f$ is 
\begin{equation}\label{Hessiano}
\,({\rm Hess}\,f)  \,(U,V)=-g\big(R(U,X)X,V\big) +
g\big(A_X(U),A_X(V)\big),
\end{equation}
for all $U,V\in \mathfrak{X}(M)$. Notice that
 $X$ is Killing if and only if  $A_X$ is
skew-adjoint with respect to $g$, and, in this case, the gradient of $f$ 
 is $A_X(X)$.
\vspace{0.5mm}

\noindent {\em Proof of Theorem  \ref{tecnica}.}
(A)  To check that $X$ is Killing, 
note that, as it is homothetic, $\lambda \,g(X,X)=Xg(X,X)=2X(f)$. The latter vanishes at $p_0$ by hypothesis, and $\lambda=0$ follows because $X_{p_0}$ is timelike.

\vspace{0.5mm}

The spacelike subspace
$X^\bot_{p_{_{_0}}}$ of $T_{p_{_{_0}}}M$ is now $A_X$-invariant and has
${\rm dim}\,X^{\bot}_{p_{_{_0}}}=2n-1$. Therefore, the induced operator
$A'_X$ of $X^\bot_{p_{_{_0}}}$ must have at least an eigenvalue $\mu\in \R$, and $\mu=0$ 
because $A'_X$ is skew-adjoint. Hence, a corresponding eigenvector $v \in
X^\bot_{p_{_{_0}}}$, ($v \neq 0$), satisfies $A_X(v)=0$.

\vspace{0.5mm}

Now, claim that $({\rm Hess} f)_{p_{_{_0}}}$ must be positive
semi-definite (because $p_{_{_0}}$ is assumed to be a local minimum). Consequently, formula  (\ref{Hessiano}) yields 
$\mathcal{K}(\pi) \geq 0$, where  $\pi = {\rm Span} \{v,
X_{p_{_{_0}}}\}$.  

\vspace{0.5mm}

(B)  As  $X$ is  homothetic, $\lambda \,g(X,Y)=g(\nabla_XX, Y) +Y(f)$ for all $Y$.  The last term vanishes at $p_0$ by hypothesis and so $A_{X}(X_{p_0})=-\lambda X_{p_0}$, that is, $X_{p_0}$ is an eigenvector of $A_{X}$ with eigenvalue $\mu=-\lambda$.

\vspace{0.5mm}

The subspace
$X^\bot_{p_{_{_0}}}$ of $T_{p_{_{_0}}}M$ has ${\rm dim}\,X^{\bot}_{p_{_{_0}}}=2n$ and is degenerate. Since $X^{\bot}_{p_{_{_0}}}$ is $A_X$-invariant, then $A_X$ can be induced to an operator  $\widehat A_{X}$ on the $(2n-1)$-dimensional quotient vector space $\widehat X^{\bot}_{p_{_{_0}}}:=X^{\bot}_{p_{_{_0}}}/$ Span$\{X_{p_{_{_0}}}\}$. Note that the metric $\hat g_{p_{_{_0}}}$ induced from the restriction $g_{p_{_{_0}}}|_{X^{\bot}_{p_{_{_0}}}}$ to $\widehat X^{\bot}_{p_{_{_0}}}$ is positive definite and $\widehat A_{X}$ is skew-adjoint with respect to $\hat g_{p_{_{_0}}}$. Reasoning as in the case (A), $\widehat A_{X}$ must have at least an eigenvector $\hat v:= v+$ Span $\{X_{p_{_{_0}}}\}$ with eigenvalue  $\mu=0$. 

Summing up, there exists $v \in
X^\bot_{p_{_{_0}}}$, with $\{v,X_{p_{_{_0}}}\}$ independent, satisfying $A_X(v)\in $ Span$\{X_{p_{_{_0}}}\}$  and, 
as $({\rm Hess} f)_{p_{_{_0}}}$ must be positive
semi-definite,
formula  (\ref{Hessiano}) yields 
$\mathcal{K}_{X_{p_0}}(\pi) \leq 0$ for  $\pi = {\rm Span} \{v,
X_{p_{_{_0}}}\}$.  $\hfill \square$

In order to prove Theorem \ref{t_3}, recall first that the Grassmannian of all the timelike tangent planes of a Lorentzian manifold $M$ is arc-connected. In fact, the timelike planes at the tangent space of a point $p\in M$ are arc-connected because both, the set of timelike directions and the set of spacelike directions are arc-connected; this includes the case $m=2$, as we are interested in the (non-oriented) directions spanned by the vectors rather than in the vectors theirselves. Moreover, any timelike plane $\pi$ tangent at  $p$ can be connected by a piecewise smooth arc to some
timelike plane tangent at a different point $q\in M$, just by propagating parallely a basis of $\pi$. Arc-connectedness  is trivially extended to the Grassmannian of causal planes and the Grassmannian of lightlike planes   (the latter, necessarily, in dimension $m\geq 3$),  as well as the  Grassmannian of timelike (resp. causal, lightlike) planes containing a  timelike (resp. causal, lightlike) $X\in \mathfrak{X}(M)$. Notice that if $X$ is causal, the Grassmannian of causal planes containing $X$ only has timelike planes whenever $X$ is timelike. 

\vspace{0.5mm}

\noindent {\em Proof of Theorem \ref{t_3}.} As the manifold is compact, $f$ must attain its minimum value at some $p$ and its maximum value at some $q$. For the case (A),  Theorem \ref{tecnica}(A) implies the existence of a timelike tangent plane $\pi_p$ such that $\mathcal{K}(\pi_p) \geq  0$. Reasoning with formula (\ref{Hessiano}) as above (or using the much more general result \cite[Th. 3.1(i)]{RS5})
yields easily that $\mathcal{K}(\pi_q) \leq  0$ for all the timelike tangent planes containing $X_q$. Thus, the result follows from the claimed arc-connectedness of the Grassmannian of timelike planes containing $X$. 

The case (B) is analogous. Indeed, now the maximum $q$ satisfies $f(q)=0$ and either reasoning with formula (\ref{Hessiano})
or using \cite[Th. 3.1]{RS5}
one has $\mathcal{K}_{X_q}(\pi_q) \geq  0$. If $f(p)<0$ then Theorem 1(A) implies $\mathcal{K}(\pi_q) \geq  0$. Thus,  the connectedness of the  Grassmannian of causal planes $X$ in addition of 
the property of continuity of timelike and null sectional curvatures (established in 
Remark \ref{segundo}(a)) imply that either \eqref{e_th} or \eqref{e_th2} must fail. If $f(p)=0$, $X$ is everywhere lightlike and Theorem 1(B) (in addition to  connectedness and continuity along the  Grassmannian of  lightlike planes on $X$) makes the job. 
$\hfill\square$

\begin{rem} {\rm Instead using \cite[Th. 3.1(i)]{RS5} to prove Theorem \ref{t_3}(A) in the case of positive timelike sectional curvature, without any restriction on the dimension, this can be also deduced from \cite[Th. 2]{HKS} where it is shown that if a timelike complete Lorentzian manifold has positive timelike sectional curvature, then it admits no timelike Killing vector field. }
\end{rem}

\section{Further discussions and examples }
\label{s3}

\subsection{The conformal case} 
As a difference with other  techniques mentioned above, which are
extendible to the conformal case \cite{RS2, RS4}, Theorem
\ref{tecnica} cannot be extended directly to the case $X$ conformal.
\begin{counter-exe}\label{conformal_extension} {\rm  Let us consider the Lorentzian metric on $\R^2$ conformally equivalent to the
$2$-dimensional Lorentz-Minkowski one, $g=e^{2u(x,y)}(dx^2-dy^2)$.
Clearly, $X=\partial_y$, ($\partial_y$ denotes $\partial /\partial y$), is a timelike conformal vector field.
Moreover, the Gauss curvature of $g$ satisfies
$$\mathcal{K}(x,y)=-e^{-2u(x,y)}\left(
\partial ^2 u/\partial x^2 -\partial^2 u/\partial y^2 \right)(x,y),$$
at any $(x,y)\in \R^2$. In order to get the required counter-example, it is enough to choose
$u(x,y) = -(x^2+2y^2)$ and to observe that the function $g(\partial_y,\partial_y)=-e^{-2(x^2+2y^2)}$ attains a minimum at $p_0=(0,0)$. However, $\mathcal{K}(x,y)<0$ for any $(x,y)\in \R^2$. }
\end{counter-exe}

\vspace{.5mm}

Despite this, one can use previous techniques to obtain some information on the conformal case. 
For example, assume there is  a conformal vector field $X$ on an even dimensional Lorentzian manifold $(M,g)$, with $\mathcal{L}_Xg=\sigma g$, and that there exists a critical point ${p_{_0}}\in M$ of the function $f=\frac{1}{2}\,g(X,X)$ such that $X_{p_{_0}}$ is timelike. 
Then, $\sigma(p_{_0})=0$ holds and the spacelike subspace $X^{\perp}_{p_{_0}}$  of $T_{p_{_0}}M$ is $A_X$-invariant. The induced operator  $A'_X$ on $X^{\perp}_{p_{_0}}$ is then skew-adjoint. Consequently, there exists $v\in X^{\perp}_{p_{_0}}$, $g(v,v)=1$, such that $A'_X(v)=0$. On the other hand, we have
$$({\rm Hess}\,f)_{p_{_0}}  \,(v,v)=-g\big(R(v,X_{p_{_0}})X_{p_{_0}},v\big) - \frac{1}{2}\,X_{p_{_0}}(\sigma).$$
Therefore, if $f$ attains a local minimum at $p_{_0}$, then we get $$\mathcal{K}(\pi)\geq \frac{1}{2} X_{p_{_0}}(\sigma)\big(-g(X_{p_{_0}},X_{p_{_0}})\big)^{-1},$$
 where $\pi={\rm Span}\{v,X_{p_{_0}}\}$. 
 Summing up: 
a sufficient extra assumption to ensure $\mathcal{K}(\pi)\geq 0$  (so extending Theorem \ref{tecnica}) is $ X_{p_{_0}}(\sigma) \geq 0$.

\begin{rem}\label{ultimo1} {\rm   Consistently with this result,  we have $ X_{p_{_0}}(\sigma) < 0$ in  Counter-example \ref{conformal_extension} (namely, $\sigma =2\, \partial u/\partial y$, $\sigma(0,0)=0$ and $\partial \sigma /\partial y(0,0)=-4$).
}
\end{rem}

\subsection{Other examples and applications}

The necessity of even $m$ in Theorem \ref{t_3}(A) is stressed by the examples  in \cite{KR, Ka} (with constant  negative sectional curvature)
pointed out in the Introduction. Next, we construct a simple  different  class.  

\begin{eje}\label{ultimo2} {\rm 
Let
$\S^{2n+1}=\big\{z:=(z_1,...,z_{n+1})\in \C^{n+1}\, :\, \sum_{j=1}^{n+1}|z_j|^2=1\big\}$ be the $(2n+1)$-dimensional unit sphere and let $X\in \mathfrak{X}(\S^{2n+1})$ be given by $X_z=iz$, for any $z\in \S^{2n+1}$.
Denote by $g^o$ the canonical Riemannian metric of $\S^{2n+1}$. The vector field $X$ is Killing for $g^o$ and satisfies $g^o(X,X)=1$, hence $\nabla^{o}_XX=0$ where $\nabla^o$ is the Levi-Civita connection of $g^o$. Making use of $g^o$ and $X$, define on $\S^{2n+1}$ the following Lorentzian metric $g=g^o-2\,\omega\otimes \omega$, where $\omega$ is the one form $g^o$-equivalent to $X$. Clearly, $X$ is also Killing for $g$ and satisfies $g(X,X)=-1$, hence $\nabla_XX=0$, where $\nabla$ is the Levi-Civita connection of $g$. On the other hand, $\nabla_UV=\nabla^o_UV-2\,\omega(U)\nabla^o_VX-2\,\omega(V)\nabla^o_UX$, for all $U,V\in \mathfrak{X}(\S^{2n+1})$,  \cite[(4.2)]{GPR}. 
 
In $(\S^{2n+1},g)$ any timelike tangent plane containing $X$ has negative sectional curvature. In fact, we have $\nabla_VX=-\nabla^o_VX$ for all $V$. Moreover, if $g(V,X)=0$ is also fulfilled, then $\nabla^o_VX=iV$ (the distinction between the covariant derivative of $\S^{2n+1}$ and the one of $\mathbb{C}^{n+1}$ is irrelevant here) and, consequently,  $\nabla_VX=-iV$. Therefore, from formula (\ref{Hessiano}) we get $\mathcal{K}(\pi)=-1$, for any tangent plane $\pi$ which contains $X$. This is in agreement to the fact that $(\S^{2n+1},g)$ is spatially isotropic with respect to $X$, i.e., at every $p \in \S^{2n+1}$, for any two unit vectors $v_1, v_2\in X_p^{\perp}$, there exists a $g$-isometry $\phi$ of $\S^{2n+1}$ such that $\phi(p)=p$, $d\phi_p(X_p)=X_p$ and $d\phi_p(v_1)=v_2$. Moreover, the unitary group $U(n+1)$ acts transitively by $g$-isometries on $\S^{2n+1}$,  \cite[Prop. 4.2]{GPR}, but $(\S^{2n+1},g)$ does not have constant sectional curvature; indeed, it is  geodesically complete (by any of the results in [M,Ka,RS2]) as well as simply connected, but not a model space-form.

For the sake of completeness,  we point out that the natural projection $\pi$ from 
$\S^{2n+1}$ onto the complex projective space $\mathbb{C}P^n$ is a semi-Riemannian submersion from $(\S^{2n+1},g)$ onto $(\mathbb{C}P^n,g_{FS})$, where $g_{FS}$ is the classical Fubini-Study metric of constant holomorphic sectional curvature $4$. Consequently, O'Neill's formula \cite[Th. 7.47]{O'N} may be called to relate the sectional curvature in $(\S^{2n+1},g)$ of a horizontal tangent plante $P=\mathrm{Span}\{u,v\}$, $g(u,u)=g(v,v)=1, g(u,v)=0$, with the sectional curvature of its projected plane $d\pi(P)$ in $(\mathbb{C}P^n,g_{FS})$. Namely, 
$\mathcal{K}_g(P)=\mathcal{K}_{g_{FS}}(d\pi(P))+3\,g(iu,v)^2$, showing also that $\mathcal{K}_g$ is non constant (recall that $\mathcal{K}_{g_{FS}}$ is non constant and satisfies $1 \leq \mathcal{K}_{g_{FS}}\leq 4$). Note that if $g$ is replaced by $g^o$, O'Neill's formula looks as follows $\mathcal{K}_{g{^o}}(P)=\mathcal{K}_{g_{FS}}(d\pi(P))-3\,g^o(iu,v)^2$, in accordance that the fact that $(\S^{2n+1},g^o)$ has constant sectional curvature $1$.
} \end{eje}

The possibilities  of changing the sign of $f=\frac{1}{2}\,g(X,X)$ for a  Killing vector field $X$ and the consistency of positive or negative sectional curvature at the local minima or maxima of $f$ can be checked by the following very simple $2$-dimensional example and $3$-dimensional variants. 

\begin{eje}\label{movido1} {\rm Let  $M=\R^2$ endowed
with the Lorentzian metric  $g= dx\otimes dy + dy\otimes dx + 2 f(x)\, dy^2$, where $f = \frac{1}{2}\,g(\partial_y,\partial_y)$ is any smooth function independent of $y$, and the Gauss curvature at each $(x_0,y_0)$ is   $-f''(x_0)$ (for this result and background, see \cite{Sa1}, formula (3.1) and Remark 9.2 (2)). If   $f$ is negative at some point $x_0$, then the Killing vector field $\partial_y$ becomes timelike at $x_0$. Consistently if,  additionally, $f$ attains at $x_0$ a local minimum, Theorem \ref{tecnica} implies that the Gauss curvature of $g$ will be non negative at $(x_0,y), y\in \R$.  Notice that, if $f$ is $1$-periodic, then  a compact Lorentzian manifold is obtained just by taking the quotient torus $\R^2/\Z^2$. The possibilities for the null sectional curvature can also be checked with variants of this example such as $M=\R^3$, $\bar g= dx\otimes dy + dy\otimes dx + 2 f(x)\, dy^2+ h(x,z)dz^2$ (see also Remark \ref{r_4}(b)).
} 
\end{eje}

Finally, we show the consistency of our bounds in the case of warped products,  which provides an illustration for the relativistic case involving Schwarzschid spacetime and its extension in Kruskal one.

\begin{eje} {\em The curvature  of a warped product $B\times_\rho F:=(B\times F, g:= g_B+ \rho^2 g_F)$, $\rho\in C^\infty(B)$, $\rho>0$, is well-understood and manageable,  as O'Neill's formulas for submersions simplify \cite[Ch. 7]{O'N}.  Assume that the base $(B,g_B)$ is Riemannian and the fiber $(F,g_F)$ is Lorentzian. Any timelike Killing vector field $X^F$ on the fiber $F$ lifts to a  Killing one $X$ on $M=B\times F$. Taking into account that $2f(x,y):= g(X,X)_{(x,y)}=\rho(x)^2 g_F(X^F_y,X^F_y)$ for all $(x,y)\in M$, $X$ is also timelike and  Theorem \ref{tecnica}(A) can be applied.
 
In particular, this  holds when  $(F,g_F)=(\R, -dt^2)$ and $X^F$ is the natural vector field $d/dt$ which lifts to $X=\partial_t$. These are the relativistic {\em static spacetimes}, which are physically and  geometrically significant (see for example \cite{Sa4}). The inequality  Ric$(v,v)\geq 0$ for any timelike tangent vector $v$, that is called the {\em timelike convergence condition} (TCC),  means that gravity, on average, attracts, while the equality (Ricci flatness) means that spacetime is vacuum therein \cite[p. 123]{SW}. 
If  $f=\frac{1}{2}\,g(\partial_t,\partial_t)$ attains a local maximum  at $p$, then
$\mathcal{K}(\pi)\leq 0$ for all timelike tangent plane $\pi$ containing $X_p$, \cite[Th. 3.1(i)]{RS5}. This fact implies Ric$(X_p,X_p)\geq 0$, in agreement with the TCC.
On the contrary, the existence of a local minimum  of $f$ at some $p$ yields, according to Theorem \ref{tecnica}, a timelike tangent plane $\pi$
containing $X_p$ with the sign of $\mathcal{K}(\pi)$ in the wrong direction for TCC. Indeed,   O'Neill formulas show that, under TCC, $K(\pi)=0$ for all the timelike planes $\pi \subset T_pM$ containing $X_p$.

Noticeably, classical Schwarzschild exterior spacetime, \cite[Def. 13.2]{O'N}, 
is  vacuum and static. It  
 satisfies
that, for each point $q$,  there is a timelike tangent plane $\pi_q$ containing $X_q$ with
$K(\pi_p)\neq 0$; thus, $f$ cannot attain a local minimum. Indeed,  $2f=2(m/r) -1$, where $r \in (2m,\infty)$ is a suitable radial coordinate of the base and the {\em mass}  $m>0$ is a constant.
Nevertheless, $f$ admits the infimum $-1$, which is ``attained'' when $r\rightarrow \infty$ (and $K(\pi_{q(r)})\rightarrow 0$). 
It is also worth pointing out that Schwarzschild exterior spacetime can be regarded as an open subset of the {\em Kruskaal spacetime}. This is not static, but the Killing vector $\partial_t$ can be extended so that it becomes lightlike in the region that would correspond to $r=2m$.
} \end{eje}

\end{document}